\documentclass[a4paper, 10pt, comsoc]{IEEEconf}

\usepackage{amsmath}
\usepackage{amssymb}
\usepackage[english]{babel}
\usepackage[utf8x]{inputenc}
\usepackage{longtable}
\usepackage{enumerate}
\usepackage{enumitem}
   \setlist{nosep}
\usepackage{url}
\usepackage{pdflscape}
\usepackage{graphicx}   
\usepackage[table]{xcolor}
\usepackage{tikz}
\usepackage{mathtools}
\usepackage{arydshln}
\usepackage{subfig}
\usepackage{keyval}
\usepackage{hhline}
\usepackage{etoolbox}
\usepackage[toc,page]{appendix}
\usepackage{comment}
\usepackage{amsthm}

\newtheorem{thm}{Theorem}

\newtheorem{prop}{Proposition}
\newtheorem{rem}{Remark}
\newtheorem*{pf*}{Proof}

\DeclareMathOperator{\Col}{col} 
\makeatletter
\newsavebox{\@brx}
\newcommand{\llangle}[1][]{\savebox{\@brx}{\(\m@th{#1\langle}\)}%
  \mathopen{\copy\@brx\kern-0.5\wd\@brx\usebox{\@brx}}}
\newcommand{\rrangle}[1][]{\savebox{\@brx}{\(\m@th{#1\rangle}\)}%
  \mathclose{\copy\@brx\kern-0.5\wd\@brx\usebox{\@brx}}}
\makeatother    \newtheorem{defn}{Definition}

\renewenvironment{abstract}{%
\hfill\begin{minipage}{0.95\textwidth}
\rule{\textwidth}{1pt}}
{\par\noindent\rule{\textwidth}{1pt}
\vspace{4pt}
\end{minipage}}

\begin{document}


\title{Structure Preserving Discretization of 1D Nonlinear Port-Hamiltonian Distributed Parameter Systems} 


\author{B.C. van Huijgevoort, S. Weiland, H.J. Zwart \\ \small Eindhoven University of Technology \\ \small E-mail: b.c.v.huijgevoort@tue.nl}




\twocolumn[
\begin{@twocolumnfalse} 
\maketitle

\begin{abstract}   
\textbf{Abstract}              
This paper contributes with a new formal method of spatial discretization of a class of nonlinear distributed parameter systems that allow a port-Hamiltonian representation over a one dimensional manifold. A specific finite dimensional port-Hamiltonian element is defined that enables a structure preserving discretization of the infinite dimensional model that inherits the Dirac structure, the underlying energy balance and matches the Hamiltonian function on any, possibly nonuniform mesh of the spatial geometry.
\end{abstract}
\end{@twocolumnfalse}
]
\pagestyle{plain}


\section{Introduction}
\label{sec:intro}
First principle mathematical models provide accurate descriptions of the behavior of physical phenomena, but their complexity is often prohibitive to perform analysis or to determine closed-form analytic solutions. In these cases, methods from scientific computing provide the adequate tools to convert these models into algebraic structures that allow for numerical simulations. The approximate nature of these methods is widely accepted, but the approximate solutions may no longer comply with the first principle laws of the underlying mathematical model. It is for this reason that the quest for structure preserving discretization methods has received considerable research attention. 
Work on symplectic and geometric integrators \cite{Hairer2006,Feng2010,Kotyczka2019} conserve symplecticity, many dedicated 3D meshing techniques
\cite{Chew1997,Teng2001,Alliez2001} preserve structure and symmetry in spatial and temporal discretization, while a wealth of dedicated post-processing tools in finite element, finite volume and finite difference methods \cite{Krist1997,Ascher1998,Johnson2012}
aim to certify compliance of numerical outcomes with the underlying physical laws. 
Distributed parameter port-Hamiltonian systems are of specific interest in this context. Starting with the early work of \cite{Courant1990,Dorfman1987,Dorfman1993,Zwart2012}, a wide range of research directions on the modeling, simulation, control and discretization  of port-Hamiltonian systems have emerged.
See, e.g. \cite{Duindam2014,VanDerSchaft2002,rashad2020twenty,Yoshimura2006,Villegas2007,Cervera2007,Farle2013} and references therein.

The discretization in space and time of distributed parameter port-Hamiltonian systems involves fundamental questions on energy conservation, power balances, modularity in system composition and the preservation of Hamiltonian properties in the discretization process.  Specific for this model class, discrete differential geometries were introduced in \cite{Talasila2004,Talasila2006} to define discrete Dirac structures and discrete port-Hamiltonian dynamics on discrete manifolds. Collocation based methods and simplicial Dirac structures are employed in \cite{Seslija2011,Seslija2014,Kotyczka2019,Kotyczka2018a,Kotyczka2018b} to derive a discrete exterior calculus in which numerical integration preserves Dirac structures, dynamics and constitutive equations. The authors of \cite{Serhani2019,Macchelli2014} consider discretizations of damped and controlled systems. Partitioned finite element methods were considered in \cite{Cardoso2018} and involve integration by parts on a subset of equations.
Methods in  \cite{Golo2004,Macchelli2007,Voss2009,Voss2011b,Voss2011s,Bassi2006,Harshit2020} focus on the spatial discretization of boundary controlled systems and construct particular finite element modules from the direct approximation of the differential forms. 

This paper is in line with the latter approaches and addresses the direct spatial discretization of nonlinear distributed parameter Hamiltonian systems (with a decomposable Hamiltonian) controlled via ports at the boundary of a 1D spatial geometry.  We propose
a novel structure preserving discretization method that (i) preserves power balances through a discrete Dirac structure, (ii) introduces generic finite-elements for non-uniform grids that locally express the dynamics of the system through a Galerkin type of projection, (iii) provides a consistent approximation of the Hamiltonian function and (iv) is valid for nonlinear distributed and boundary controlled port-Hamiltonian systems.

We show that the Stokes-Dirac structure underlying the distributed port Hamiltonian system can be projected on a finite dimensional differential form over the 1D spatial manifold so as to establish a constant and discrete Dirac structure that is fully separated from the dynamics and the definition of the Hamiltonian of the underlying distributed parameter model. The technical novelty to establish this, lies in the introduction of a fictitious point inside the mesh elements that partition the spatial geometry. This augments the degrees of freedom and remedies a number of limitations concerning compatibility conditions and choices of input-output variables \cite{Golo2004}, direct feedthrough terms that cause oscillations due to instantaneous power flows through the discretized network \cite{Voss2009,Voss2011b,Voss2011s}, non-sparsity or non feasibility to match nonlinear Hamiltonian functions \cite{Harshit2020}. As such, the proposed discretization scheme is numerically efficient, flexible towards choices of input, output and boundary variables, fully scalable to large networks and physically relevant on its preservation of passivity and energy distribution properties.

The paper is organized as follows. A concise problem definition is given in Section~\ref{sec:pf}. State space formulations of lumped- and distributed-parameter port-Hamiltonian systems are given in Section~\ref{sec:pHsystems}. The discretization method and the interconnection structure are described in Section~\ref{sec:discretization}. Section~\ref{sec:illustration} illustrates the method on an application of a lossless transmission line. Conclusion and recommendations are given in Section~\ref{sec:conclusions}. 

\section{Problem formulation}
\label{sec:pf}
In abstract form, a Dirac structure is a subspace 
$\mathcal{D}$ of a Cartesian product $\mathcal{F}\times\mathcal{E}$ where $\mathcal{F}$ is a linear space and $\mathcal{E}=\mathcal{F}^*$ its dual
with respect to the pairing $\langle \cdot \mid \cdot \rangle: \mathcal{E}\times\mathcal{F}\rightarrow\mathbb{R}$. The spaces $\mathcal{E}$ and $\mathcal{F}$ are usually referred to as the spaces of \emph{efforts} and \emph{flows}, respectively, while $\langle e \mid f\rangle$, the evaluation of the linear functional $e\in\mathcal{E}$ on $f\in\mathcal{F}$, expresses the \emph{power}. The Euclidean product $\mathcal{F}\times\mathcal{E}$ defines the space of \emph{port variables} and carries an indefinite symmetric bi-linear pairing $\llangle \cdot,\cdot\rrangle$ defined as
\begin{align}
\label{eq:bilineardef}
\begin{split}
&\llangle (f_1,e_1),(f_2,e_2) \rrangle := \langle e_1 \mid f_2 \rangle + \langle e_2 \mid f_1 \rangle, \\ 
&(f_i,e_i)\in\mathcal{F}\times\mathcal{E}.
\end{split}
\end{align} 
\begin{defn}
\label{def:Dirac}
A (constant) Dirac structure on $\mathcal{F}\times \mathcal{E}$ is a subspace $\mathcal{D}\subset \mathcal{F}\times\mathcal{E}$ such that
$\mathcal{D} = \mathcal{D}^{\perp}$,
where $\mathcal{D}^{\perp}:=\{d'\in\mathcal{D}\mid \llangle d,d'\rrangle =0 \textrm{ for all } d\in\mathcal{D} \}$ denotes the orthogonal of $\mathcal{D}$ with respect to the bi-linear pairing $\llangle \cdot, \cdot \rrangle$ defined in \eqref{eq:bilineardef}.
\end{defn}
In particular, any $(f,e)$ belonging to a Dirac structure satisfies $0=\llangle(f,e),(f,e)\rrangle=2\langle e \mid f\rangle$ showing that Dirac structures have power conserving properties over their port variables.

A Hamiltonian is a function $H:\mathcal{X} \rightarrow\mathbb{R}$ defined on a differential manifold  $\mathcal{X}$ in which $H(x)$ represents the energy in a storage element $x\in\mathcal{X}$. A port-Hamiltonian system has partitioned flows $\mathcal{F}=\mathcal{F}_s\times \mathcal{F}_e$ and efforts $\mathcal{E}=\mathcal{F}_s^*\times\mathcal{F}_e^*=
\mathcal{E}_s\times\mathcal{E}_e$ on which the indefinite symmetric bi-linear pairing \eqref{eq:bilineardef} extends to 
\begin{align*}
  \label{eq:bilinearext}
  \begin{split}
  &\llangle (f_1,e_1),(f_2,e_2)\rrangle 
  = \\
  &\langle e^s_1 \mid f^s_2\rangle +
  \langle e^s_2 \mid f^s_1\rangle + 
  \langle e^e_1 \mid f^e_2\rangle + 
  \langle e^e_2 \mid f^e_1\rangle.
   \end{split}
   \end{align*}
Here, $\mathcal{F}_s$ is the flow space of tangent vectors of (any) $x\in\mathcal{X}$, $\mathcal{F}_e$ the external flow space and  $(f^s,e^s)$ and $(f^e,e^e)$ refer to the \emph{storage} and \emph{external} port variables of the system, respectively. The dynamics of a port-Hamiltonian system is fully defined by $(\mathcal{D},H,\mathcal{X},\mathcal{F}_e)$ in the sense 
that the storage port variables $(f^s,e^s) = (-\dot{x},\delta_x H(x))$ where $\delta_x H(x)$ is the gradient or the variational derivative of the Hamiltonian $H$ at the point $x$, (depending on whether $\mathcal{X}$ is finite or infinite dimensional, resp.), and $(f^e,e^e)$ represents the port variables through which the system interacts with its environment. Specifically, the dynamics is defined by  
\begin{equation}
   \label{eq:dyndirac}
   \left(-\dot{x}(t),\delta_x H(x(t)),f^e(t),e^e(t) \right) \in \mathcal{D}.
\end{equation}
for all time $t\geq 0$. Whenever time courses of \eqref{eq:dyndirac} are well defined and continuously differentiable, the Dirac structure $\mathcal{D}$ imposes the power balance
\begin{equation}
    \label{eq:conservation}
    \frac{\text{d}H}{\text{d}t} = 
    \langle \delta_xH(x(t)) \mid \dot{x}(t)\rangle = \langle e^e(t) \mid f^e(t) \rangle,\quad t\geq 0
\end{equation}
which expresses that the storage rate of a port-Hamiltonian system equals the power supplied through its external port variables. It is emphasized that, although the Dirac structure is linear, the dynamics \eqref{eq:dyndirac} becomes nonlinear whenever $\delta_xH(x)$ is nonlinear in $x$. 

A port-Hamiltonian distributed parameter system has an infinite dimensional flow space $\mathcal{F}=\mathcal{F}_s\times\mathcal{F}_e$, where $\mathcal{F}_s$ is a space of functions on a compact set $Z\subset\mathbb{R}^n$, representing the spatial domain, extended with  functions on the boundary $\partial Z$ of $Z$. It is the purpose of this paper to approximate a port-Hamiltonian distributed-parameter system  $\Sigma=(\mathcal{D},H,\mathcal{X}, \mathcal{F}_e)$ by a port-Hamiltonian lumped-parameter system $\Sigma_N=(\mathcal{D}_N,H_N,\mathcal{X}_{N}, \mathcal{F}_{e,N})$ that is, in fact, an aggregation of $N$ components. Ideally, this approximate system meets a number of requirements:
\begin{itemize}[leftmargin=*]
    \item \textbf{Requirement 1.} 
    $\mathcal{D}_N$ defines a Dirac structure on a finite dimensional space $\mathcal{F}_N\times\mathcal{E}_N$ which is the canonical projection $\pi_N = (\pi_{N\mathcal{F}},\pi_{N\mathcal{E}})$ of $\mathcal{F}\times\mathcal{E}$ onto 
    $\mathcal{F}_N\times\mathcal{E}_N\subset\mathcal{F}
    \times \mathcal{E}$, where $\mathcal{E}_N$ is the dual of $\mathcal{F}_N$ with respect to the pairing 
    $\langle e_N\mid f_N\rangle_N := 
    \langle e_N\mid f_N\rangle$ for any $e_N\in\mathcal{E}_N\subset\mathcal{E}$ and $f_N\in\mathcal{F}_N\subset\mathcal{F}$ and 
    \begin{equation}
        \label{eq:diracincl}
         \pi_N\mathcal{D}^{\perp}=\pi_N\mathcal{D}\quad\subseteq\quad\mathcal{D}_N = \mathcal{D}_N^{\perp}
    \end{equation}
    where $\mathcal{D}_N^{\perp}$ is the orthogonal of $\mathcal{D}_N$ with respect to the bi-linear pairing
    \[
       \llangle (f_{N}^1,e_{N}^1), (f_{N}^2,e_{N}^2)\rrangle_N:=
       \langle e_{N}^1\mid f_{N}^2\rangle_N +\langle e_{N}^2\mid f_{N}^1\rangle_N.
    \]
    In words, $\mathcal{F}_N$ inherits the duality pairing from $\mathcal{F}$ while the inclusion \eqref{eq:diracincl} implies that for any \mbox{$(f,e)\in\mathcal{D}$}, its projection $(f_N,e_N):=\pi_N(f,e)$ is power conserving in the aggregated system in the sense that $\langle f_N \!\mid\! e_N\rangle_N \!=\! 0$. 
    \item \textbf{Requirement 2.} 
    The Hamiltonian $H_N$ is defined on a finite dimensional state space $\mathcal{X}_{N}\subset\mathcal{X}$ and achieves
    \[
         H_N(x)=H(x) \quad\text{for all }x\in\mathcal{X}_{N} 
    \]
    i.e., the total storage $H$ of the distributed parameter system when restricted to the finite dimensional state space $\mathcal{X}_{N}$ matches the total storage in the aggregated system.
\end{itemize}
The combined properties therefore imply that the power balance \eqref{eq:conservation} of the distributed system is preserved 
in the aggregated system $\Sigma_N$ in the sense that 
$\dot{H}_N(t) =\langle f_N^e(t),e_N^e(t)\rangle$ for all $t\geq 0$.  The problem that is solved in this paper amounts to explicitly constructing, for any $N>0$, an aggregated port-Hamiltonian lumped-parameter system $\Sigma_N$ that approximates the  port-Hamiltonian distributed-parameter system $\Sigma$ while meeting these two requirements. 

\section{Port-Hamiltonian systems in state space coordinates}
\label{sec:pHsystems}

\subsection{Lumped parameter port-Hamiltonian systems}
\label{ssec:pHlumped}
Finite-dimensional port-Hamiltonian systems are generally given in input-state-output form by the equations
\begin{equation}
\begin{bmatrix} \dot{x}\\ y\end{bmatrix} 
= \begin{bmatrix} A & B \\ C & D \end{bmatrix} \begin{bmatrix} \delta_xH(x)\\
u\end{bmatrix}
\label{eq:pHlumped}
\end{equation} 
where $x(t)\in\mathcal{X}=\mathbb{R}^{n_x}$ is the state variable and $H:\mathcal{X}\rightarrow \mathbb{R}$ is the Hamiltonian function. The system interacts with its environment through the input-output pair $(u, y)$, where $\dim(u)=\dim(y)=n_e$. 
The real-valued matrices $(A,B,C,D)$ are of compatible dimensions and satisfy
\begin{equation}
\label{eq:Rdef}
R:=\begin{bmatrix}-A-A^{\top} & C^{\top}-B \\ C-B^{\top} & D+D^{\top}
\end{bmatrix} \succeq 0.
\end{equation}
Here, $R$ represents the losses in the system \cite{Duindam2014} while \eqref{eq:Rdef} implies that the dissipation inequality
\begin{equation}
\frac{\textrm{d}}{\textrm{d}t}H(x(t)) \leq y(t)^{\top}u(t)
\label{eq:DIE}
\end{equation}
holds for all time $t\geq 0$ and all system trajectories $(u,x,y)$ compatible with \eqref{eq:pHlumped}. If \eqref{eq:DIE} is an equality, the system is said to be \emph{conservative}, which is the case if and only if $R=0$.

Following the abstract setting of Section~\ref{sec:pf}, \eqref{eq:pHlumped} is fully defined by the Hamiltonian $H$, the storage port variables $(f^s,e^s):=(-\dot{x},\frac{\partial H}{\partial x})$ and the external port variables  $(f^e,e^e):=(u,y)$. Indeed, the algebraic relation
\begin{equation}
\begin{bmatrix}f^s \\ e^e\end{bmatrix} =
\begin{bmatrix}-A & -B \\ C & D \end{bmatrix}
\begin{bmatrix} e^s \\ f^e\end{bmatrix} 
\label{eq:DiracLumped}
\end{equation}
satisfies the \emph{power balance}
\begin{equation}
 e^{s \top}f^s + e^{e \top}f^e = \frac{1}{2}\begin{pmatrix}e^s \\ f^e \end{pmatrix}^{\top} R\begin{pmatrix}e^s \\ f^e \end{pmatrix}.
 \label{eq:Plumped}
\end{equation}
with $R$ defined in \eqref{eq:Rdef}.
In particular, with $\mathcal{F}:=\mathcal{F}_s\times\mathcal{F}_e=\mathbb{R}^{n_x+n_e}$, the space $\mathcal{E}=\mathbb{R}^{n_x+n_e}$ is dual to $\mathcal{F}$ with respect to the pairing 
\begin{equation}
  \langle e \mid f\rangle =e^{\top}f=
  \langle e^s\mid f^s\rangle +\langle e^e\mid f^e\rangle
  \label{eq:pairing}
\end{equation}
and the following result on conservative systems is an immediate consequence of \eqref{eq:Plumped}. 
\begin{prop}
If $R=0$, then 
\begin{align}
\mathcal{D}:=\left\{((f^s,f^e),(e^s,e^e))\in \mathcal{F}\times\mathcal{E} \mid \eqref{eq:DiracLumped} \text{ holds} \right\}.
\label{eq:Diracstruc}
\end{align}
is a finite-dimensional Dirac structure and the system \eqref{eq:pHlumped} is fully defined by $(\mathcal{D},H,\mathcal{X},\mathcal{F}_e)$ according to \eqref{eq:dyndirac}.
\label{prop:FiniteDirac}
\end{prop}
Hence, if $R=0$, \eqref{eq:pHlumped} defines a finite dimensional Dirac structure \eqref{eq:Diracstruc} while, conversely, \eqref{eq:Diracstruc}, together with $H$, fully defines \eqref{eq:pHlumped}. 

\begin{rem}
If $R\succeq 0$ is non-zero, the port-Hamiltonian systems \eqref{eq:pHlumped} can also be cast in the setting of Dirac structures. Specifically, factorize $R=G^{\top}R_0G$ in which $R_0\succ 0$ and $G=\left[\begin{smallmatrix}G_x \\ G_e \end{smallmatrix}
\right] \in \mathbb{R}^{(n_x+n_e)\times n_r}$ with $n_r=\mathrm{rank}R$. Define $\mathcal{D}$ as the graph of the skew symmetric operator
\[
   \begin{bsmallmatrix} 
   -\frac{1}{2}(A-A^{\top}) & -B & -G_x \\ B^{\top} & \frac{1}{2}(D-D^{\top}) & -G_e \\ 
   G_x^{\top} & G_e^{\top} &0 \end{bsmallmatrix}
\]
mapping $(e^s,f^e,e^r)\mapsto(f^s,e^e,f^r)$.  
Then $\mathcal{D}$ is a Dirac structure by adding the resistive port $\langle f^r,e^r\rangle$ to the pairing in \eqref{eq:pairing} and \eqref{eq:pHlumped} is represented by 
\[
   (-\dot{x},\delta_xH(x(t)),u(t),y(t),f^r(t),e^r(t))\in\mathcal{D}
\]
where $f^r=-R_0e^r$.
\end{rem}

\subsection{Distributed parameter port-Hamiltonian systems}
\label{ssec:pHdistributed}
Distributed parameter port-Hamiltonian systems are described by partial differential equations whose solutions evolve over time and an $n$-dimensional bounded spatial domain $Z\subset \mathbb{R}^n$. Generally, these systems can be written in the form
\cite{Duindam2014}
\begin{equation}
\begin{bmatrix}
\frac{\partial p}{\partial t} \\
\frac{\partial q}{\partial t}
\end{bmatrix} =
A \begin{bmatrix}
\delta_p H \\
\delta_q H
\end{bmatrix},\quad
\begin{bmatrix}
u \\
y
\end{bmatrix} = 
\begin{bmatrix}
B \\
C
\end{bmatrix}\begin{bmatrix}
\delta_p H |_{\partial Z} \\
\delta_q H |_{\partial Z}
\end{bmatrix},
\label{eq:pHdist}
\end{equation} 
where the energy variables $p(z,t)$ and $q(z,t)$ depend on a spatial coordinate $z\in Z$ and time $t$. Furthermore, $A$ is an operator whose domain is a dense subset of a (Hilbert) space of functions on $Z$ with the property that $A+A^*\preceq 0$ in the sense that $\langle x,Ax\rangle + \langle A^*x,x\rangle \leq 0$ for all compactly supported elements $x$ in the domain of $A$, and $\delta_p$ and $\delta_q$ denote the variational derivatives with respect to the energy variables $p$ and $q$, respectively. The variables $u(z,t)$ and $y(z,t)$ are the input and output  
on the boundary $z\in\partial Z$ for time $t\geq 0$. 

With $\Lambda^k(Z)$, $k=0,\ldots, n$ denoting the space of exterior $k$-forms on the manifold $Z$, the Hamiltonian $H:\Lambda^{n_p}(Z)\times\Lambda^{n_q}(Z) \rightarrow 
\mathbb{R}$ describes the total energy in the system and is given by $H(p,q) = \int_Z \mathcal{H}(p,q)$ \cite{Duindam2014,VanDerSchaft2002}
where $\mathcal{H}(p,q)$ is the energy density. Throughout, it is assumed that $\mathcal{H}$ can be decomposed according to 
\begin{equation}
   \label{eq:cHdef}
   \mathcal{H}(p,q) = \mathcal{H}_p(p) + \mathcal{H}_q(q),
\end{equation}
where $\mathcal{H}_p:\Lambda^{n_p}(Z)\rightarrow\Lambda^{n}(Z)$,
$\mathcal{H}_q:\Lambda^{n_q}(Z)\rightarrow\Lambda^{n}(Z)$ and
$n_p+n_q = n+1$. 

Following the setting of 
Section~\ref{sec:pf}, the efforts and flows associated with \eqref{eq:pHdist} are defined by
\begin{subequations}
\label{eq:efDefinf}
\begin{alignat}{2}
f^p &:= -\frac{\partial p}{\partial t},&\quad
e^p &:= \delta_p H = \frac{\partial \mathcal{H}_p(p)}{\partial p},
\label{eq:pDefinf}\\
f^q &:= -\frac{\partial q}{\partial t},&\quad
e^q &:= \delta_q H = \frac{\partial \mathcal{H}_q(q)}{\partial q}.
\label{eq:qDefinf}
\end{alignat}
\end{subequations} 
Since $f^p\in\Lambda^{n_p}(Z)$, $f^q\in\Lambda^{n_q}(Z)$, $e^p\in\Lambda^{n-n_p}(Z)$, $e^q\in\Lambda^{n-n_q}(Z)$, this motivates defining the linear spaces of flows and efforts as
\begin{subequations}
\label{eq:spaceFE}
\begin{align}
\label{eq:spaceF}
\mathcal{F}_{p,q} &:= \Lambda^{n_p}(Z) \times \Lambda^{n_q}(Z) \times \Lambda^{n-n_p}(\partial Z) \\
\label{eq:spaceE}
\mathcal{E}_{p,q} &:= \Lambda^{n-n_p}(Z) \times \Lambda^{n-n_q}(Z) \times \Lambda^{n-n_q}(\partial Z). 
\end{align}
\end{subequations} 
In particular, $\mathcal{E}_{p,q}$ is the dual of $\mathcal{F}_{p,q}$ with respect to the non-degenerate pairing 
\begin{align}
   \label{eq:pHDirac-pairing}
   &\langle(e^p,e^q,e^b) \mid (f^p,f^q,f^b)\rangle :=\\
   &\int_Z[e^p \wedge f^p + e^q \wedge f^q] + 
   \int_{\partial Z} [e^b \wedge f^b] \nonumber
\end{align}
where $\wedge$ is the usual wedge product of differential forms. Here, the first term in the right-hand side of \eqref{eq:pHDirac-pairing} represents the power in the distributed storage port and the second term the power delivered over the boundary $\partial Z$ of the domain \cite{Duindam2014,Golo2004,VanDerSchaft2002}. 
In turn, \eqref{eq:pHDirac-pairing} equips $\mathcal{F}_{p,q}\times\mathcal{E}_{p,q}$ with the symmetric bi-linear form \eqref{eq:bilineardef}.  Partition $\mathcal{F}_{p,q}=\mathcal{F}_s\times\mathcal{F}_b$ with
$\mathcal{F}_s:=\Lambda^{n_p}(Z)\times\Lambda^{n_q}(Z)$ the storage flow space and $\mathcal{F}_b=\Lambda^{n-n_p}(\partial Z)$ the external flow space. Then, $(f^b, e^b)\in\mathcal{F}_b\times \mathcal{F}_b^*$ denote the port variables that represent the interaction of the system at the boundary $\partial Z$ and coincide with  $(u,y)$ in \eqref{eq:pHdist}.
Using these definitions, the system in \eqref{eq:pHdist} can be written as 
\begin{align}
\begin{bmatrix}
-\frac{\partial p}{\partial t} \\
-\frac{\partial q}{\partial t} 
\end{bmatrix} = 
\begin{bmatrix}
\sigma \star & (-1)^r\text{d}\\
\text{d} & 0
\end{bmatrix} \begin{bmatrix}
\delta_p H(p,q) \\
\delta_q H(p,q)
\end{bmatrix}, \nonumber\\
\begin{bmatrix}
f^b \\
e^b
\end{bmatrix} = 
\begin{bmatrix}
1 & 0 \\
0 & -(-1)^{n-n_q} 
\end{bmatrix} \begin{bmatrix}
\delta_pH |_{\partial Z} \\
\delta_qH |_{\partial Z}
\end{bmatrix}
\label{eq:pHdistr}
\end{align}
where  $r=n_pn_q+1$, $\star:\Lambda^k(Z)\rightarrow\Lambda^{n-k}(Z)$ is the Hodge star operator\footnote{The Hodge star operator satisfies $\omega_1 \wedge \star \omega_2 = \langle \omega_1\mid \omega_2 \rangle \omega$.} and $\text{d}$ the exterior derivative. Losses in the system are represented by $\sigma$ which is a scalar if the dissipation is homogeneous in all directions, it is a diagonal matrix if the dissipation is inhomogeneous and isotropic and it is a non-symmetric matrix if the dissipation is an-isotropic. The Hamiltonian $H$ in \eqref{eq:pHdistr} satisfies the Stokes-Cartan expression for $t\geq 0$
    \begin{align*}
     \frac{\textrm{d} H}{\textrm{d}t} &= 
     \int_Z \left[\delta_pH \wedge \frac{\partial p}{\partial t}+ \delta_qH\wedge \frac{\partial q}{\partial t}\right] \\
     &= -\!\!\int_{\partial Z}[(-1)^{n-n_q}\delta_qH \wedge \delta_pH]  -\!\! \int_Z [\delta_p H \wedge \sigma \star \delta_pH],
     \end{align*}
or, stated in terms of the effort and flow variables \eqref{eq:efDefinf},
\begin{align}
   \int_Z \! [e^p \!\wedge\! f^p \!+\! e^q \!\wedge\! f^q ] \!+\!
   \int_{\partial Z}\![e^b \!\wedge\! f^b] \!=\! \int_Z \![e^p \!\wedge\! \sigma\star e^p] 
\label{eq:Pdistr}
\end{align}
which is the power balance for \eqref{eq:pHdistr}, generalizing \eqref{eq:Plumped} to the infinite-dimensional case.

\begin{prop}
If $\sigma=0$, then 
\begin{align}
\mathcal{D} :=  & \left\{(f^p,f^q,f^b),(e^p,e^q,e^b) \in \mathcal{F}_{p,q} \times \mathcal{E}_{p,q} \mid \right.\nonumber \\
 & 
\begin{bmatrix} f^p \\ f^q \end{bmatrix} = 
\begin{bmatrix}
0 & (-1)^r\text{d}\\ \text{d} & 0 \end{bmatrix} 
\begin{bmatrix} e^p \\ e^q \end{bmatrix},\nonumber\\
& \left.
\begin{bmatrix} f^b \\ e^b \end{bmatrix} = 
\begin{bmatrix} 1 & 0 \\ 0 & -(-1)^{n-n_q} \end{bmatrix} 
\begin{bmatrix} e^p |_{\partial Z} \\ e^q |_{\partial Z} \end{bmatrix}
\right\}.
\label{eq:Diracdist}
\end{align} 
is a Dirac structure on $\mathcal{F}_{pq}\times\mathcal{E}_{pq}$ and the system \eqref{eq:pHdistr} is fully defined by
$(\mathcal{D},H,\mathcal{X},\mathcal{F}_e)$ according to \eqref{eq:dyndirac}.
\label{prop:InfiniteDirac}
\end{prop} 
\begin{pf*}
    This Proposition is proven in \cite{VanDerSchaft2002}. 
\end{pf*}

Thus, if $\sigma=0$ the conservative port-Hamiltonian system \eqref{eq:pHdistr} satisfies the Dirac structure \eqref{eq:Diracdist}. Conversely, 
\eqref{eq:Diracdist} together with $H$ fully defines  \eqref{eq:pHdistr}.

\begin{rem}
If $\sigma\neq 0$ the Dirac structure \eqref{eq:Diracdist} can be extended so as to represent the non-conservative port-Hamiltonian system \eqref{eq:pHdistr} also in in terms of an extended Dirac structure together with the Hamiltonian function $H$. See, e.g., Chapter 6 in \cite{Villegas2007}.
\end{rem}

\section{Discretization in 1D}
\label{sec:discretization}
This paper focuses on one-dimensional manifolds. That is, the spatial configuration space $Z$ is assumed to be a bounded and closed subset of $\mathbb{R}$ with a zero-dimensional boundary $\partial Z$. For one-dimensional manifolds the pair $(n_p,n_q)=(1,1)$. It is assumed that a conservative distributed-parameter port-Hamiltonian system $\Sigma=(\mathcal{D},H,\mathcal{X},\mathcal{F}_b)$ is given, as described in Section~\ref{sec:pHsystems}, 
where the spaces \eqref{eq:spaceFE} equal
\begin{subequations}
\label{eq:spaceFE1D}
\begin{alignat}{2}
\label{eq:spaceF1D}
&\mathcal{F}_{p,q} &:= \Lambda^1(Z) \times \Lambda^1(Z) \times \Lambda^0(\partial Z)&=\mathcal{F}_s\times\mathcal{F}_b\\
\label{eq:spaceE1D}
&\mathcal{E}_{p,q} &:= \Lambda^{0}(Z) \times \Lambda^{0}(Z) \times \Lambda^{0}(\partial Z)&=\mathcal{E}_s\times\mathcal{E}_b.
\end{alignat}
\end{subequations} 
The idea is to construct an aggregated lumped-parameter port-Hamiltonian model \newline   $\cramped{\Sigma_N\!\!=\!\!(\mathcal{D}_N,\!H_N,\! \mathcal{X}_{N},\! \mathcal{F}_{e,N})}$, meeting the requirements stated in Section~\ref{sec:pf}, as the interconnection of $N$ component systems, each of which is a lumped-parameter port-Hamiltonian model
     $\cramped{ \Sigma_{i}:=(\mathcal{D}_{ab,i},H_{ab,i}, \mathcal{X}_{ab,i}, \mathcal{F}_{ab,e,i})},$ for $i\!\!=\!\!1,\ldots, N$
that approximates $\Sigma$ locally at the domain $Z_{a_ib_i}$, the $i$th element in a simplicial decomposition or partition of $Z$ into $N$ elements. The structure of a mesh element assumes the form $\cramped{Z_{a_ib_i}:=[a_i,b_i]}$, with $\cramped{a_i<b_i}$. The state space of the aggregated model is $\cramped{\mathcal{X}_{N}=\Pi_{i=1}^N \mathcal{X}_{ab,i}}$ while the Hamiltonian $\cramped{H_N(x_1,\ldots, x_N) = \sum_{i=1}^NH_{ab,i}(x_i)}$ represents the total energy in the aggregated system with $H_{ab,i}$ the Hamiltonian of the component system $\Sigma_i$.

\subsection{Dirac structure of the component systems} 
To construct $\Sigma_i$, we consider a generic element $Z_{ab}$ in the partition of $Z$ and define a projection $\pi_{ab}$ of $\mathcal{F}_{p,q}\times\mathcal{E}_{p,q}$ onto 
$\mathcal{F}_{ab}\times\mathcal{E}_{ab}$ where
\begin{subequations}
\begin{alignat*}{2}
&\mathcal{F}_{ab} &:= \Lambda^1(Z_{ab}) \times \Lambda^1(Z_{ab}) \times \Lambda^0(\partial Z_{ab})\\
&\mathcal{E}_{ab} &:= \Lambda^{0}(Z_{ab}) \times \Lambda^{0}(Z_{ab}) \times \Lambda^{0}(\partial Z_{ab}),
\end{alignat*}
\end{subequations} 
the localized flow and effort spaces \eqref{eq:spaceFE1D}. Then 
$\mathcal{E}_{ab}$ is the dual of $\mathcal{F}_{ab}$
with respect to the duality pairing
\begin{align}
   \label{eq:pHDirac-ab-pairing}
   &\langle(e^p,e^q,e^b) \mid (f^p,f^q,f^b)\rangle_{ab} :=\\
   &\int_{Z_{ab}}[e^p \wedge f^p + e^q \wedge f^q] + 
   \int_{\partial Z_{ab}} [e^b \wedge f^b] \nonumber
\end{align} 
and we wish to establish the inclusion \eqref{eq:diracincl} in Requirement 1, but localized to $Z_{ab}$, i.e.
\begin{equation}
\label{eq:diracincl1D}
  \pi_{ab}\mathcal{D}^{\perp}=\pi_{ab}\mathcal{D}\ \subset\ \mathcal{D}_{ab}=\mathcal{D}_{ab}^{\perp}.
\end{equation}
To do this, let $\cramped{Z_{ab}\!=\![a,b]}$ be a 
mesh element and suppose that $m$ is such that $\cramped{a<m<b}$. For an arbitrary flow $\cramped{(f^p,f^q,f^b)\in \mathcal{F}_{p,q}}$ the projection  $\cramped{\pi_{ab,\mathcal{F}}=(\pi_{ab,\mathcal{F}}^p,\pi_{ab,\mathcal{F}}^q, \pi_{ab,\mathcal{F}}^b)}$ results in the storage flow variables
\begin{subequations}
  \label{eq:fApprox1D}
  \begin{alignat}{2}
     \label{eq:fpApprox1D}
      & f_{ab}^p(z) &:= \pi_{ab,\mathcal{F}}^pf^p(z) &= f_{am}^p \omega_{am}^p(z) +
        f_{mb}^p \omega_{mb}^p(z) \\
     \label{eq:fqApprox1D} 
      & f_{ab}^q(z) &:= \pi_{ab,\mathcal{F}}^qf^q(z) &= f_{am}^q \omega_{am}^q(z) +
        f_{mb}^q \omega_{mb}^q(z),
\end{alignat} 
\end{subequations}
with the coefficient vector
\begin{equation}
   \label{eq:fdef}
    f^s:=\Col(f_{am}^p,f_{mb}^p,f_{am}^q,f_{mb}^q)\in\mathbb{R}^4
\end{equation}
and where the one-form shape functions $\omega_{s}^p$ and $\omega_{s}^q$ indexed by $s\in\{am,mb\}$ satisfy
\begin{equation}
\int_{Z_{s_1}} \omega_{s_2} = \begin{cases}
1 \text{ for } s_1 = s_2 \\
0 \text{ for } s_1 \neq s_2.
\end{cases} 
\label{eq:1form1}
\end{equation} 
Thus, the shape functions decouple the different flows on the internal line segments $Z_{am}$ and $Z_{bm}$. Similarly, for $(e^p,e^q,e^b)\in \mathcal{E}_{p,q}$, the projection \mbox{$\cramped{\pi_{ab,\mathcal{E}}=(\pi_{ab,\mathcal{E}}^p,\pi_{ab,\mathcal{E}}^q, \pi_{ab,\mathcal{E}}^b)}$} produces the effort storage variables 
\begin{subequations}
  \label{eq:eApprox1D}
  \begin{alignat}{2} 
     \label{eq:epApprox1D}
      &e_{ab}^p\!(z)\!\!:=\! \pi_{ab,\mathcal{E}}^pe^p\!(z)\!\! =\!\!  
      e_{a}^p \omega_{a}^p(z)\! + \!
         e_m^p \omega_m^p\!(z) \!+ \!e_{b}^p \omega_{b}^p(z) \\
     \label{eq:eqApprox1D}
      &e_{ab}^q\!(z)\!\!:=\!\pi_{ab,\mathcal{E}}^qe^q\!(z) \!\! =\!\! e_{a}^q \omega_{a}^q\!(z) \!+ \!e_m^q \omega_m^q\!(z)\! +\! e_{b}^q \omega_{b}^q(z), 
   \end{alignat}
\end{subequations} 
with the coefficient vector
\begin{equation}
  \label{eq:vdef}
   v := \Col(e_a^p,  e_m^p, e_b^p, e_{a}^q, e_m^q, e_{b}^q) \in \mathbb{R}^6
\end{equation}
and the zero-form shape functions $\omega_{s}^p$, $\omega_{s}^q$ indexed by $s\in\{a,m,b\}$ satisfy
\begin{align}
\label{eq:0form}
 &\omega_{a}(a) =1,\quad\omega_{a}(b)=0,\quad\omega_{b}(a)=0, \quad\omega_b(b)=1 \\
 &\omega_{m}(a) =0,\quad\omega_{m}(m)=1, \quad\omega_m(b)=0. \nonumber
\end{align}
Given the shape functions satisfying \eqref{eq:1form1} and \eqref{eq:0form}, define the matrices
\begin{align}
M_1 & := 
\begin{bmatrix}
\int_{Z_{am}} \sigma \star \omega_a^p & \int_{Z_{am}} \sigma \star \omega_m^p & \int_{Z_{am}} \sigma \star \omega_b^p \\
\int_{Z_{mb}} \sigma \star \omega_a^p & \int_{Z_{mb}} \sigma \star \omega_m^p & \int_{Z_{mb}} \sigma \star \omega_b^p 
\end{bmatrix},
\label{eq:M11D}\\
M_2&:= 
\begin{bmatrix}
\int_{\partial Z_{am}} \omega_{a}^q & \int_{\partial Z_{am}}  \omega_{m}^q & \int_{\partial Z_{am}}  \omega_{b}^q \\
\int_{\partial Z_{mb}} \omega_{a}^q & \int_{\partial Z_{mb}}  \omega_{m}^q & \int_{\partial Z_{mb}}  \omega_{b}^q 
\end{bmatrix},
\label{eq:M21D}\\
M_3&:= 
\begin{bmatrix}
\int_{\partial Z_{am}} \omega_{a}^p & \int_{\partial Z_{am}}  \omega_{m}^p & \int_{\partial Z_{am}} \omega_{b}^p \\
\int_{\partial Z_{mb}} \omega_{a}^p & \int_{\partial Z_{mb}}  \omega_{m}^p & \int_{\partial Z_{mb}} \omega_{b}^p 
\end{bmatrix}.
\label{eq:M31D}
\end{align}
Then the following result relates the coefficient vectors $f^s$ and $v$ defined in \eqref{eq:fdef} and \eqref{eq:vdef}.
\begin{prop}
\label{Prop1}
The projections \eqref{eq:fApprox1D} and \eqref{eq:eApprox1D} satisfy 
\[
   \int_{Z_{ab}}\begin{bmatrix}f_{ab}^p(z) \\ f_{ab}^q(z)\end{bmatrix}=\int_{Z_{ab}}\begin{bmatrix} \sigma\star& \text{d}\\ \text{d} & 0\end{bmatrix} \begin{bmatrix}e_{ab}^p(z) \\ e_{ab}^q(z)\end{bmatrix} \textrm{ if and only if }
\]
\begin{equation}
 \label{eq:Mfdef}
   f^s=M_fv \quad \text{ where }
   M_f = 
  \begin{bmatrix}
      M_1 & M_2 \\ M_3 & 0
  \end{bmatrix}.
\end{equation}
\end{prop}

\begin{pf*}
Substitute the expansions \eqref{eq:fApprox1D} and \eqref{eq:eApprox1D} into
the residual expression
\[
r(z)=\begin{bmatrix}r^p(z)\\ r^q(z)\end{bmatrix} =\begin{bmatrix}f_{ab}^p(z) \\ f_{ab}^q(z)\end{bmatrix} -\begin{bmatrix} \sigma\star& \textrm{d}\\ \textrm{d} & 0\end{bmatrix} \begin{bmatrix}e_{ab}^p(z) \\ e_{ab}^q(z)\end{bmatrix}.
\]
This gives the 1-forms $r^p(z)\! =\!f_{am}^p \omega_{am}^p(z)\! +\! f_{mb}^p \omega_{mb}^p(z)
     \! -\! \sigma\! \star\! [e_{a}^p \omega_{a}^p(z) 
     \!+\!e_m^p \omega_m^p(z) +\! e_{b}^p \omega_{b}^p(z)]
     \!-\!\text{d}[e_{a}^q \omega_{a}^q(z)\! +\! e_m^q \omega_m^q(z)\! +\! e_{b}^q   \omega_{b}^q(z)]$ and $r^q(z)\! =\!f_{am}^q \omega_{am}^q(z)\! +\! f_{mb}^q \omega_{mb}^q(z) 
-\text{d}[e_{a}^p \omega_{a}^p(z)\! +\! e_m^p \omega_m^p(z) +\! e_{b}^p \omega_{b}^p(z)].$ \newline 
(If): By using Stokes' theorem and the condition \eqref{eq:1form1} on the shape functions, integration over the line segments $Z_{am}$ and $Z_{mb}$ yields $\cramped{\int_{Z_{ab}} r^p=0}$ with $M_1$ and $M_2$ given in \eqref{eq:M11D} and \eqref{eq:M21D}. Similarly integration of $r^q$ over the line segments $Z_{am}$ and $Z_{mb}$ while using Stokes' theorem and \eqref{eq:1form1} yields  $\cramped{\int_{Z_{ab}}r^q=0}$ with $M_3$ given in \eqref{eq:M31D}. \newline
(Only if): Conversely, $\cramped{\int_{Z_{ab}}r=\int_{Z_{am}}r+\int_{Z_{mb}}r =0}$ implies, by \eqref{eq:1form1} and Stokes' theorem, that $f^s=M_fv$ with $M_f$ defined in \eqref{eq:Mfdef}.
\end{pf*}
In the following, we define a power balance for the expansions \eqref{eq:fApprox1D} and \eqref{eq:eApprox1D} that matches \eqref{eq:pHDirac-ab-pairing} and the inclusion \eqref{eq:diracincl1D}. For this, first observe that, by
\eqref{eq:Diracdist} and \eqref{eq:pHDirac-ab-pairing},
\begin{align*}
&\langle (e_{ab}^p,e_{ab}^q,e_{ab}^b) \mid (f_{ab}^p,f_{ab}^q,f_{ab}^b) \rangle_{ab}=\\ &
\int_{Z_{ab}} [e_{ab}^p \wedge f_{ab}^p + e_{ab}^q \wedge f_{ab}^q] - \int_{\partial Z_{ab}} [e_{ab}^q \mid_{\partial Z_{ab}} \wedge e_{ab}^p\mid_{\partial Z_{ab}}].
\end{align*}
In the conservative case where $\sigma=0$, the latter expression vanishes by \eqref{eq:Pdistr}. 
Define the matrices
\begin{align}
M_4 &:=
\int_{Z_{ab}} \begin{bmatrix}
\omega_a^p \wedge \omega_{am}^p & \omega_m^p \wedge \omega_{am}^p & \omega_b^p \wedge \omega_{am}^p \\
\omega_a^p \wedge \omega_{mb}^p & \omega_m^p \wedge \omega_{mb}^p & \omega_b^p \wedge \omega_{mb}^p 
\end{bmatrix}
\label{eq:M41D} \\
 M_5 &:=
\int_{Z_{ab}} \begin{bmatrix}
\omega_a^q \wedge \omega_{am}^q & \omega_m^q \wedge \omega_{am}^q & \omega_b^q \wedge \omega_{am}^q\\
 \omega_a^q \wedge \omega_{mb}^q & \omega_m^q \wedge \omega_{mb}^q& \omega_b^q \wedge \omega_{mb}^q
\end{bmatrix}
\label{eq:M51D} \\
M_6 &:= \int_{\partial Z_{ab}}\!\! \begin{bmatrix}
\omega_a^q\!\mid_{\partial Z} \\
\omega_m^q\!\mid_{\partial Z} \\
\omega_b^q\!\mid_{\partial Z}
\end{bmatrix} \!\!\wedge\! 
\begin{bmatrix}
\omega_a^p\mid_{\partial Z} \
\omega_m^p\mid_{\partial Z} \
\omega_b^p\mid_{\partial Z}
\end{bmatrix}\!.
\label{eq:M61D}
\end{align} 
We then have the following expressions for the distributed and boundary power. 
\begin{prop}
\label{prop:powerchar}
The projections \eqref{eq:fApprox1D} and \eqref{eq:eApprox1D}
with $f^s$ defined by \eqref{eq:fdef} achieve that  
\begin{enumerate}
    \item the distributed power
      \begin{equation}
        \label{eq:powerDistApprox}
        P_{ab}^s:=\int_{Z_{ab}} e_{ab}^p \wedge f_{ab} ^p + e_{ab}^q \wedge f_{ab}^q = e^{s \top}f^s
      \end{equation}
      if and only if 
      \begin{equation}
        \label{eq:Medef}
         e^s = M_ev\quad \text{ where } M_e=\begin{bmatrix}M_4 & 0 \\ 0 & M_5 \end{bmatrix}.
      \end{equation}
    \item the boundary power
      \begin{equation}
          \label{eq:powerBoundApprox}
          P_{ab}^b:= -\int_{\partial Z_{ab}} e_{ab}^q\mid_{\partial Z_{ab}} \wedge e_{ab}^p\mid_{\partial Z_{ab}} = e^{e\top}f^e 
      \end{equation}
      if $f^e=M_uv$ and $e^e=M_yv$ where
      \begin{subequations}
        \label{eq:Muydef}
        \begin{align}
         M_u&:= \begin{bmatrix}0 & 0 & 0 & -1 & 0 & 0\\
         0 & 0 & 1 & 0 & 0& 0\end{bmatrix}
         \begin{bmatrix}M_6 & 0 \\ 0 & I \end{bmatrix}  \\
         M_y&:= \begin{bmatrix}1 & 0 & 0 & 0 & 0 & 0\\
         0 & 0 & 0 & 0 & 0& -1\end{bmatrix}
         \begin{bmatrix}M_6 & 0 \\ 0 & I \end{bmatrix}.
        \end{align}
      \end{subequations} 
\end{enumerate}
\end{prop}
\begin{pf*}
(1) Substituting the expansions  \eqref{eq:fApprox1D} and \eqref{eq:eApprox1D} in the definition of $P_{ab}^s$ gives
\begin{align*}
\int_{Z_{ab}} e_{ab}^p(z)\wedge f_{ab}^p(z)  
= \begin{bmatrix} e_{am}^p \\ e_{mb}^p \end{bmatrix}^{\top}
\begin{bmatrix} f_{am}^p \\ f_{mb}^p
\end{bmatrix},\\
\int_{Z_{ab}} e_{ab}^q(z)\wedge f_{ab}^q(z)
=\begin{bmatrix}
e_{am}^q \\ e_{mb}^q
\end{bmatrix}^{\top}\begin{bmatrix}
f_{am}^q \\ f_{mb}^q
\end{bmatrix},
\end{align*}
if and only if $e^s:=\Col(e_{am}^p, e_{mb}^p, e_{am}^q, e_{mb}^q)=M_ev$ is defined through $M_4$ and $M_5$ as in \eqref{eq:Medef}. This yields \eqref{eq:powerDistApprox}.

(2) Substitute the effort expansion \eqref{eq:eApprox1D} in the definition of the boundary power $P_{ab}^b$ to obtain that 
$P_{ab}^b = - \Col(e_{a}^q, e_m^q, e_{b}^q)^{\top}
M_6 \Col(e_a^p, e_m^p, e_b^p).$
Point $m$ lies in the interior of $Z_{ab}$ which, by \eqref{eq:0form}, implies that $\omega_m\mid_{\partial Z} = 0$. Thus \eqref{eq:0form} implies that $P_{ab}^b$ does not depend on the coefficient $e_m^q$. With $f^e:=M_uv$, $e^e:=M_yv$ and \eqref{eq:Muydef} it then follows that $P_{ab}^b= v^{\top}M_y^{\top}M_uv=e^{e\top}f^e$ which is \eqref{eq:powerBoundApprox}.
\end{pf*}

With the pairing \eqref{eq:pHDirac-ab-pairing}, Proposition~\ref{prop:powerchar} therefore establishes the power balance
\[
  \langle  (e_{ab}^p,e_{ab}^q,e_{ab}^b) \mid (f_{ab}^p,f_{ab}^q,f_{ab}^b)\rangle_{ab}= e^{s\top}f^s + e^{e\top}f^e.
\]
The choice of boundary port variables $(f^e,e^e)$ as defined through \eqref{eq:Muydef} is certainly not unique in establishing item 2 of Proposition~\ref{prop:powerchar}. Using Proposition~\ref{Prop1} and Proposition~\ref{prop:powerchar}, all efforts and flows of the approximate model can now be expressed in terms of the coefficients $v$ by setting
\begin{equation}
   \label{eq:sbApprox} 
   \begin{bmatrix}
     f^s \\ e^e
   \end{bmatrix}\!=\! 
   \begin{bmatrix}
      M_f \\  M_y
   \end{bmatrix}\!v 
   =:\! E^{\top}v; 
   \quad
   \begin{bmatrix}
      e^s \\ f^e
   \end{bmatrix}\!=\! 
   \begin{bmatrix}
       M_e\\  M_u
   \end{bmatrix}\!v
  =:\! F^{\top}v
\end{equation}
where $M_f$, $M_e$, $M_y$ and $M_u$ are defined in \eqref{eq:Mfdef}, \eqref{eq:Medef} and \eqref{eq:Muydef}.
Thus, $f^s\in \mathcal{F}_s:=\mathbb{R}^4$, $e^s\in\mathcal{E}_s :=\mathbb{R}^4$, $f_e\in\mathcal{F}_e := \mathbb{R}^2$ and $e^e\in\mathcal{E}_e :=\mathbb{R}^2$ and we have the following result.
\begin{prop}
\label{prop:fdDirac}
If the system \eqref{eq:pHdistr} is conservative, i.e., if  $\sigma=0$, then there exist shape functions in the expansions
\eqref{eq:fApprox1D}, \eqref{eq:eApprox1D} satisfying \eqref{eq:1form1} and \eqref{eq:0form} such that $M_1=0$ and 
\begin{equation}
\mathcal{D}_{ab} := \textrm{Im}\begin{bmatrix}
E^{\top} \\
F^{\top}
\end{bmatrix} = \textrm{Ker} \begin{bmatrix}
F & E
\end{bmatrix}
\label{eq:ImDirac1D}
\end{equation} 
is a finite-dimensional constant Dirac structure in the sense of Definition~\ref{def:Dirac} with respect to the pairing
$\langle (e^s,e^e),(f^s,f^e)\rangle=e^{s\top}f^s + e^{e\top}f^e$.
\end{prop}

\begin{pf*}
The proof in Appendix B 
 shows that the properties of a Dirac structure, namely $FE^{\top}+EF^{\top} = 0$ and rank$\begin{bmatrix}
F \mid E
\end{bmatrix} = $\textrm{dim} $\mathcal{F}$
hold.
\end{pf*}

Proposition~\ref{prop:fdDirac} establishes the requirement \eqref{eq:diracincl1D} on the composite system. We emphasize that the proof is constructive in the choice of the shape functions.  

\subsection{Hamiltonian of the component systems}
To meet Requirement 2, define the expansions
\begin{subequations}
\label{eq:energyApprox1D}
\begin{align}
   \label{eq:pabdef}
   p_{ab}(z) &= p_{am}\omega_{am}^p(z) + p_{mb}\omega_{mb}^p(z) 
   \\
   \label{eq:qabdef}
   q_{ab}(z) &= q_{am}\omega_{am}^q(z) + q_{mb}\omega_{mb}^q(z)
\end{align}
\end{subequations}
in the coefficient vectors 
$p:=\Col(p_{am}, p_{mb})$ and $q:=\Col(q_{am},q_{mb})$
where the shape functions coincide with the ones in \eqref{eq:fApprox1D} and satisfy \eqref{eq:1form1}. 
Define the Hamiltonian function $H_{ab}:\mathcal{X}_{ab}\rightarrow\mathbb{R}$ on the state space $\mathcal{X}_{ab}:=\mathbb{R}^4$ as
\begin{equation}
   \label{eq:HamiltonianApprox1D}
   H_{ab}(p,q) := H_{ab}^p(p)+H_{ab}^q(q)
\end{equation}
where $H_{ab}^p(p) := \int_{Z_{ab}}\mathcal{H}_p(p_{ab}),
    H_{ab}^q(q) := \int_{Z_{ab}}\mathcal{H}_q(q_{ab})$
in which $p_{ab}$ and $q_{ab}$ satisfy  \eqref{eq:energyApprox1D} and the energy distribution $\mathcal{H}$ is decomposed as in \eqref{eq:cHdef}. It thus follows that
     $H_{ab}(p,q) =\int_{Z_{ab}} \mathcal{H}(p_{ab},q_{ab})$
showing that \eqref{eq:HamiltonianApprox1D} matches the Hamiltonian $H$ of the distributed system when restricted to a 4 dimensional subspace of energy variables locally defined on $Z_{ab}$.

\subsection{State space representation of component systems} \label{sec:statespace}
As stated in Section~\ref{sec:pf}, the Dirac structure $\mathcal{D}_{ab}$ and the Hamiltonian function $H_{ab}$ fully define the dynamics of the component system through the implicit relation
\[
    (-\dot{x}(t), \delta_xH_{ab}(x(t)),f^e(t),e^e(t)) \in \mathcal{D}_{ab},
    \quad t\geq 0.
\]
The matrices \eqref{eq:Muydef} are chosen such that $F$ in \eqref{eq:sbApprox} is non-singular. As a consequence, an explicit and unconstrained input-state-output representation of the component system is obtained with $f^e=u$ the input and $e^e =y$ the output:

\begin{prop}
\label{prop:pHfinCal}
The matrix $F$ in \eqref{eq:Mfdef} is non-singular and the quadruple $(A,B,C,D)$ defined as
\begin{equation}
    \begin{bmatrix}-A & -B \\ C & D \end{bmatrix}:=
    E^{\top}F^{-\top} = \begin{bmatrix}
     M_f \\ \hdashline[1pt/2pt] M_y
\end{bmatrix}\begin{bmatrix} M_e \\ \hdashline[1pt/2pt]
M_u
\end{bmatrix}^{-1}
\label{eq:pHfinal1D}
\end{equation}
satisfies \eqref{eq:Rdef}. Moreover, with $H_{ab}$ as defined in \eqref{eq:HamiltonianApprox1D} and state $x=\Col(p,q)$ the system
\[
\cramped{\Sigma_{ab}: \quad \begin{cases}
\dot{x} &= A\frac{\partial H_{ab}}{\partial x}(x) + Bu \\
y &= C\frac{\partial H_{ab}}{\partial x}(x)+Du
\end{cases}}
\]
is a $4$th order port-Hamiltonian system of the form \eqref{eq:pHlumped}.
\end{prop}

\begin{pf*}
The result is an immediate consequence of the Dirac structure from Proposition~\ref{prop:fdDirac} and is obtained by eliminating $v$ from \eqref{eq:sbApprox} through \eqref{eq:pHfinal1D}.  
\end{pf*}

\subsection{The aggregated model}

By Proposition~\ref{prop:pHfinCal}, each element $Z_{a_i,b_i}$, $\cramped{i=1,\ldots, N}$ in the partition of $Z$ defines a port-Hamiltonian system \mbox{$\Sigma_{i}:=\Sigma_{ab,i}$} defined by the Dirac structure $\mathcal{D}_i:=\mathcal{D}_{ab, i}$ set in \eqref{eq:ImDirac1D} and the Hamiltonian 
$H_{ab,i}:\mathcal{X}_{ab,i}\rightarrow\mathbb{R}$ defined in \eqref{eq:HamiltonianApprox1D}. The port-Hamiltonian structure of the aggregated model is an immediate consequence of the composition property of the Dirac structure of the component systems $\Sigma_1,\ldots, \Sigma_N$.  

\begin{figure}
    \begin{tikzpicture}[scale=0.6] 
\def\firstellipse{(0,0) ellipse (1.5 and 1)}
\def\secondellipse{(4,0) ellipse (1.5 and 1)}
\def\thirdellipse{(8,0) ellipse (1.5 and 1)}
\filldraw[fill=yellow!20]\firstellipse;
\node (A) at (0,0) {$\mathcal{D}_{i-1}$};
\filldraw[fill=yellow!20]\secondellipse;
\node (B) at (4,0) {$\mathcal{D}_i$};
\filldraw[fill=yellow!20]\thirdellipse;
\node (C) at (8,0) {$\mathcal{D}_{i+1}$};
\node at (2,-2) {$\begin{pmatrix}f^e_{i,i-1} \\ e^e_{i,i-1}\end{pmatrix}$};
\node at (6,-2) {$\begin{pmatrix} f^e_{i,i+1}\\ e^e_{i,i+1}\end{pmatrix}$};
\draw [-,thick] (4.1,1) -- (4.1,2);
\draw [-,thick] (5.5,0.1) -- (6.5,0.1);
\draw [-,thick] (3.9,1) -- (3.9,2);
\draw [-,thick] (5.5,-0.1) -- (6.5,-0.1);
\draw [-,thick] (0.1,1) -- (0.1,2);
\draw [-,thick] (-1.5,0.1) -- (-2.5,0.1);
\draw [-,thick] (1.5,0.1) -- (2.5,0.1);
\draw [-,thick] (-0.1,1) -- (-0.1,2);
\draw [-,thick] (-1.5,-0.1) -- (-2.5,-0.1);
\draw [-,thick] (1.5,-0.1) -- (2.5,-0.1);
\draw [-,thick] (8.1,1) -- (8.1,2);
\draw [-,thick] (9.5,0.1) -- (10.5,0.1);
\draw [-,thick] (7.9,1) -- (7.9,2);
\draw [-,thick] (9.5,-0.1) -- (10.5,-0.1);
\end{tikzpicture}
    \caption{Composition of Dirac structures}
    \label{fig:interconnect1D}
\end{figure}
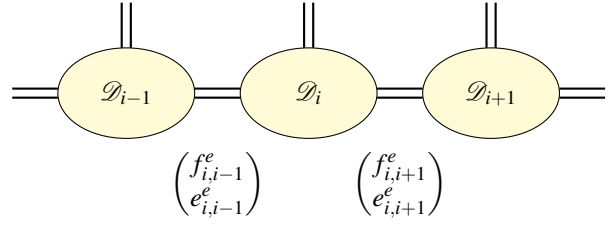

Specifically, the composition of two Dirac structures $\mathcal{D}'
\subset (\mathcal{F}'\times\mathcal{F}^e)\times (\mathcal{E}'\times\mathcal{E}^e)$ 
and $\mathcal{D}''\subset (\mathcal{F}''\times\mathcal{F}^e)\times (\mathcal{E}''\times\mathcal{E}^e)$ that share a common flow and effort space $\mathcal{F}^e$ and $\mathcal{E}^e$ is defined as
\begin{align*}
   &\mathcal{D}'\sqcap\mathcal{D}'':=
   \{ (f',f'',e',e'') \mid \exists (f^e,e^e) \text{ such that } \\
    & (f',f^e,e',e^e) \in\mathcal{D}' \text{ and }
    (f'',-f^e,e'',e^e) \in\mathcal{D}''\}.
\end{align*}
Here, the minus sign establishes that power $\langle f^e\mid e^e\rangle$ entering $\mathcal{D}'$ over the shared port is the negative of the power entering $\mathcal{D}''$, making the interconnection neutral over the shared port variables. It is easily seen that the composition
$\mathcal{D}'\sqcap\mathcal{D}''$ is again a Dirac structure with respect to the bi-linear pairing \eqref{eq:bilineardef}.

In the construction of the aggregated model, we distinguish between external and internal boundaries.  
The external boundary of $\Sigma_i$ is the (possibly empty) set $\partial Z \cap \partial Z_{a_i,b_i}$. The internal boundary of $\Sigma_i$ is $\partial Z_{a_i,b_i}\backslash\partial Z$. If $\Sigma_i$ has internal boundaries, then $Z_{a_{i-1},b_{i-1}}$ and/or $Z_{a_{i+1},b_{i+1}}$ 
can be assumed neighbouring elements to $Z_{a_i,b_i}$ in the sense that
  $a_{i-1} < b_{i-1} =a_i < b_i =a_{i+1} < b_{i+1}.$
Then $\Sigma_{i-1}$ and $\Sigma_i$ have the internal boundary point $b_{i-1}=a_i$ in common, while $\Sigma_i$ and $\Sigma_{i+1}$ have the common internal boundary point $b_i=a_{i+1}$. Let $\mathcal{D}_{i-1}$, $\mathcal{D}_{i}$ and $\mathcal{D}_{i+1}$ be the associated Dirac structures of these systems, respectively. Then the external port variables $(f_i^e,e_i^e)$ of $\mathcal{D}_i$ are partitioned with the internal boundary points $a_i$ and $b_i$ as
\begin{align*}
     f_i^e &= \begin{bmatrix}f_{i,i-1}^e \\  f_{i,i+i}^e\end{bmatrix} \in\mathbb{R}^2, \quad
     e_i^e = \begin{bmatrix}e_{i,i-1}^e \\  e_{i,i+i}^e\end{bmatrix}\in \mathbb{R}^2
\end{align*}
where $(f_{i,i-1}^e,e_{i,i-1}^e)$ and $(f_{i,i+1}^e,e_{i,i+1}^e)$
denote the port variables of $\Sigma_i$ (equivalently $\mathcal{D}_i$)
at the points $a_i$ and $b_i$ of the domain $Z_{a_ib_i}$, respectively. See Figure~\ref{fig:interconnect1D}. 
At the external boundary points, one simply matches the boundary port variables $(f^b,e^b)$ defined in \eqref{eq:pHdistr} or \eqref{eq:Diracdist}.

By composing all $N$ elements of the partitioning of $Z$, one infers a Dirac structure
\[
 \mathbb{D}_N:=\mathcal{D}_1 \sqcap \cdots \sqcap\mathcal{D}_N
\]
of the aggregate system. The Hamiltonian function $H_N:\Pi_{i=1}^N\mathcal{X}_{ab,i}\rightarrow\mathbb{R}$ is defined as the total energy
    $H_{N}(x_1,\ldots, x_N):= \sum_{i=1}^N H_{ab,i}(x_i)$.
This construction results in a system with the following properties.

\begin{thm}
  \label{thm:main}
  For every $N>0$, the aggregated system $\Sigma_N:=(\mathbb{D}_N,H_N,\mathbb{R}^{4N},\mathbb{R}^2)$ is a $4N$th order port-Hamiltonian system that satisfies Requirement 1 and Requirement 2 of Section~\ref{sec:pf}.
\end{thm}

State space representations of $\Sigma_N$ directly follow from Proposition~\ref{prop:pHfinCal} and are sparse. In particular, the kernel representation \eqref{eq:ImDirac1D} defines implicit state space models $0=F\Col(\partial_x H_{ab},u) + E\Col(-\dot{x},y)$ of the component systems which, using  \eqref{eq:sbApprox}, gives a sparsity index (ratio of zero and non-zero entries) of at least
$1- \frac{6}{6N-2}$ in the aggregate model.

\section{Simulation example}
\label{sec:illustration}
\begin{figure*}
\centering
\subfloat[Hamiltonian of this approach (blue) and of the finite element method (black)]{
\includegraphics[width=0.7\columnwidth]{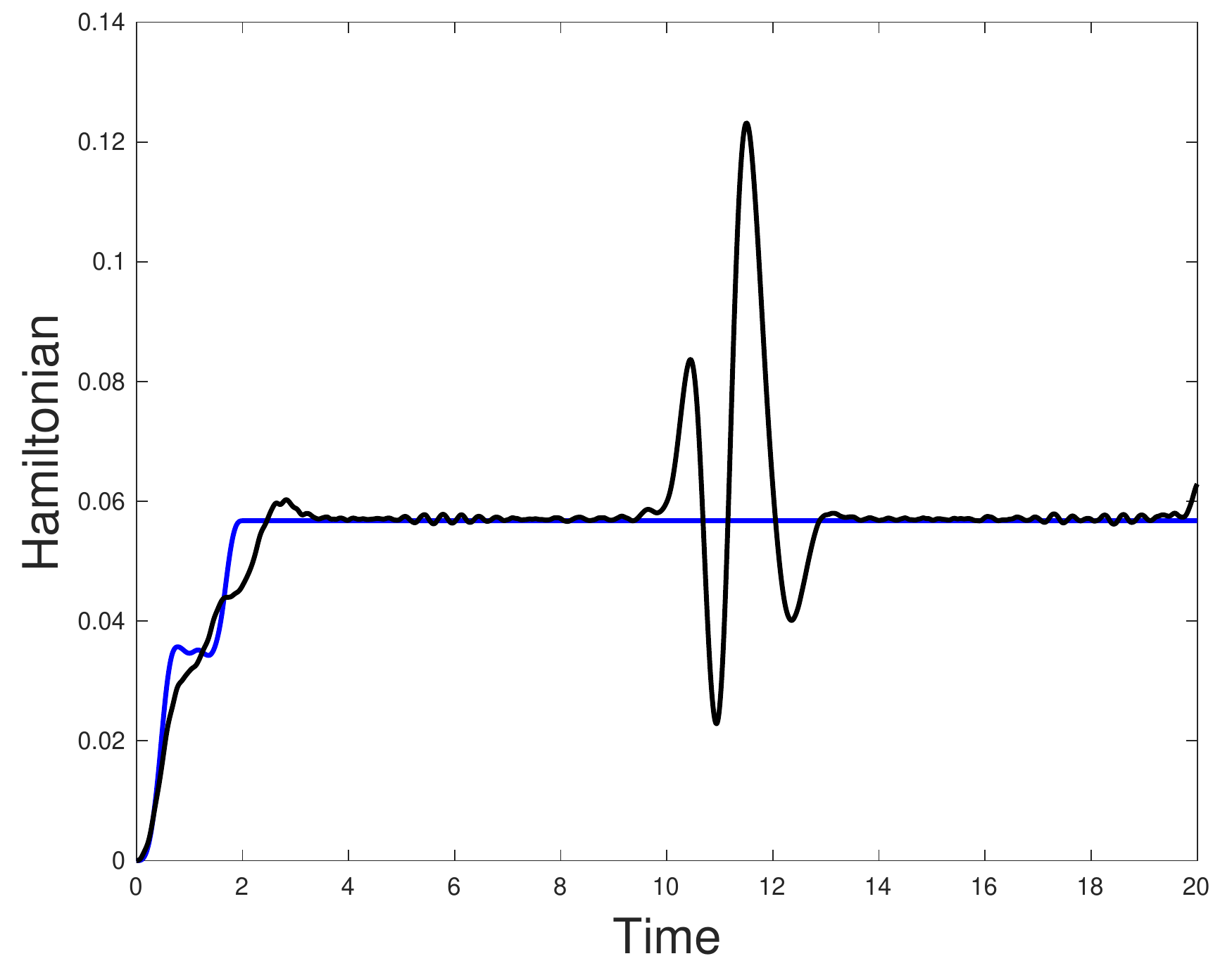}
\label{fig:ApproxH}}
\subfloat[Conservativity]{
\includegraphics[width=0.7\columnwidth]{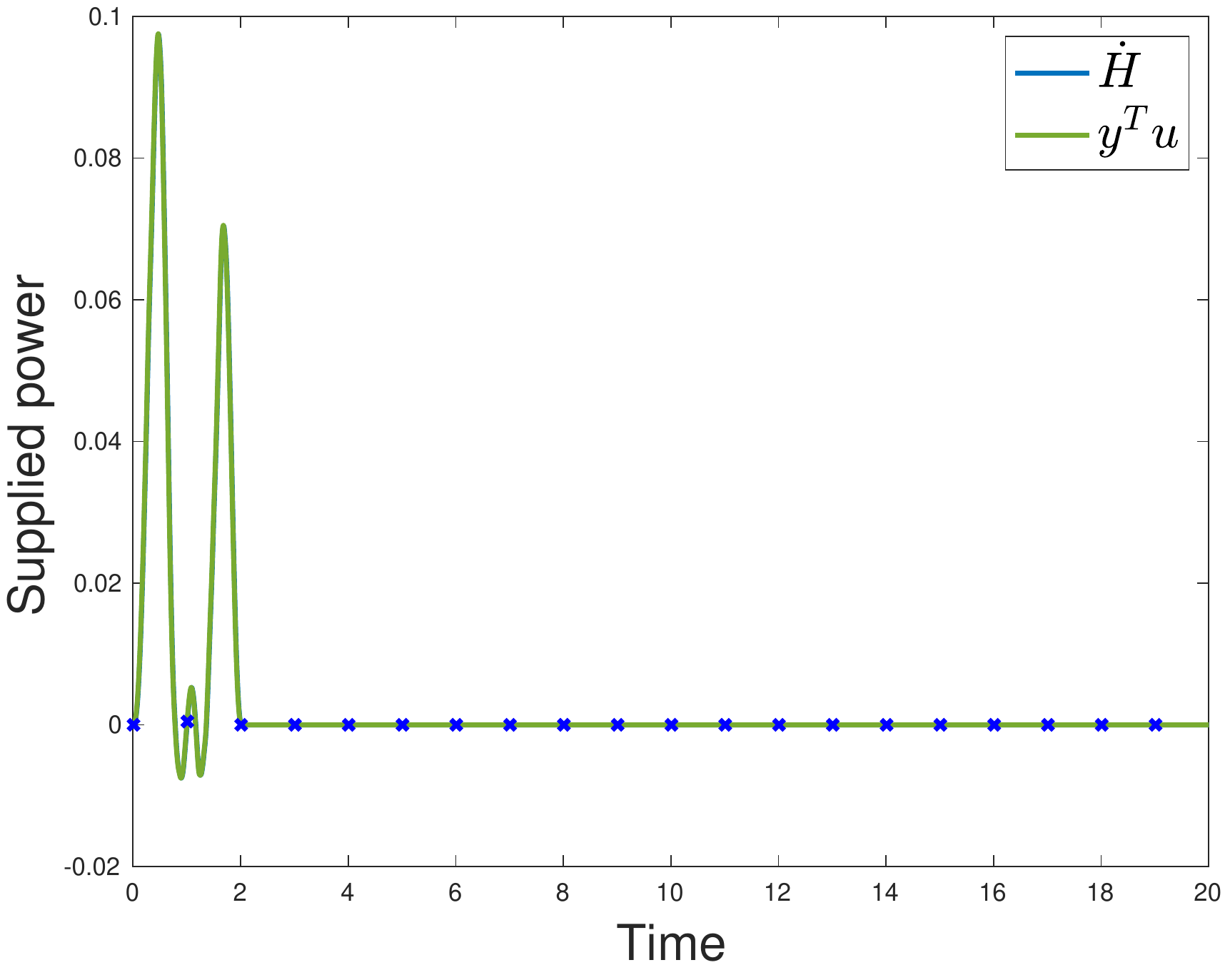}
\label{fig:Conserv}}
\caption{Preservation of the power balance}
\label{fig:H,Cons}
\end{figure*}

This section illustrates the discretization method described in this paper on the example of a lossless transmission line. The spatial-temporal voltage and current behavior of a lossless transmission line is described by 
\[
     \frac{\partial V}{\partial z} =-L\frac{\partial I}{\partial t},
     \qquad
     \frac{\partial I}{\partial z} =-C\frac{\partial V}{\partial t}
\]
where $V(t,z)$ and $I(t,z)$ denote voltage and current at time $t\geq 0 $ and position $z\in Z:=[0,\ell]$ in a transmission line of length $\ell$ with distributed capacitance and inductance $C(z)$ and $L(z)$. 
These equations can be written in port-Hamiltonian form \cite{Golo2004} by defining magnetic flux and charge distributions $p$ and $q$ as state variables 
and by introducing the Hamiltonian function
\begin{equation*}
H(p,q) = \frac{1}{2}\int_0^{\ell} \left[\frac{1}{C(z)} q \wedge \star q + \frac{1}{L(z)} p \wedge \star p\right].
\end{equation*}
This gives a lossless distributed parameter port-Hamiltonian model
\eqref{eq:pHdistr} with $(n,n_p,n_q)\!=\!(1,1,1)$, $\sigma\!=\!0$ and where the variational derivatives $e^p\!:=\!\delta_{p}H(p,q)\!=\!I$ and $e^q\!:=\!\delta_{q}H(p,q)\!=\!V$ denote current and voltage in the line and $(f^b,e^b)\!=\!(I|_{\partial Z},V|_{\partial Z})$ denote the current and voltage at either of the end-points of the line.

We let $\ell\!=\!10$ and assume uniform capacitance $C(z)\!=\!C\!=\!10^{-2}$ and uniform inductance $L(z)\!=\!L\!=\!1$. Inputs and outputs are chosen from the boundary variables $(f^b,e^b)$ according to
\[
u(t) = \begin{bmatrix}
V(t,0) \\
I(t,\ell)
\end{bmatrix},\quad
y(t) = \begin{bmatrix}
-I(t,0) \\
V(t,\ell)
\end{bmatrix}
\]
where the input $u(t)$ is  set as
\[
V(t,0) = \begin{cases}
\sin(\pi t) & \text{for } 0 \leq t \leq 2 \\
0 & \text{for } t>2
\end{cases},\quad 
I(t,\ell) = 0
\]
i.e., a one period sinusoidal voltage excitation is applied on the left-end, while the right-end of the transmission line is open.
\begin{figure}
    \centering
   \includegraphics[scale=0.5]{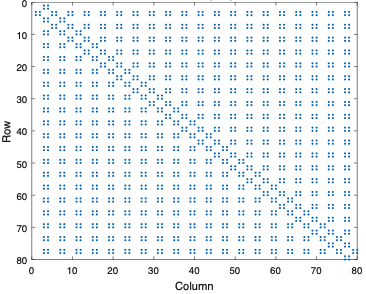}
    \caption{Sparsity of the state evolution matrix $A$; each dot represents a non-zero term}
    \label{fig:sparse}
\end{figure}
The geometry $Z=[0,\ell]$ is partitioned in $N=20$ equal line segments of length $\ell/N$. The aggregated model $\Sigma_{20}$ of Theorem~\ref{thm:main} has been constructed and has state dimension $4N=80$. Figure~\ref{fig:ApproxH} shows a time evolution of the approximated Hamiltonian. Since the input voltage has compact support in the time interval $[0,2]$, power is only delivered to the transmission line in the first 2 seconds, after which the internal energy remains constant because of the lossless characteristics of the transmission line. The approximate model clearly shows this property in Figure
~\ref{fig:ApproxH}. In order to verify whether the approximate system is indeed conservative we determine the time evolution of the power supplied to the circuit. This is computed by taking the time derivative of the Hamiltonian and by computing the pointwise product of the inputs and outputs of the system. Since the
transmission line is lossless, 
\eqref{eq:DIE} becomes an equality 
$\frac{\text{d}H}{\text{d}t} = \langle y(t) \mid u(t) \rangle$ for $t\geq 0$.
In Figure~\ref{fig:Conserv} the blue line and the red line correspond to the time derivative of the Hamiltonian and the computed power $\langle y(t) \mid u(t) \rangle$ in the aggregated model. The exact matching of the evolutions in Figure~\ref{fig:Conserv} shows that the aggregated model is indeed conservative and that the power balance is correctly preserved.
 
Even though there is a feedthrough term in the approximate model $\Sigma_{20}$, the state space representation of the interconnected system has sparse matrices.  In this simulation, the state evolution matrix $A \in \mathbb{R}^{80 \times 80}$  has a sparsity index of $\frac{4720}{80\cdot 80}=0.7375$. That is, almost three out of four entries vanish in the approximate model respresentation. Evidently, this index is beneficial for the scalability of the method when partition refinements are applied. The sparsity structure of the state evolution matrix $A$ is illustrated in Figure~\ref{fig:sparse}. 

As a comparison, the same example implemented in the Partial Differential Equation Toolbox of \texttt{Matlab}, involves a classical Finite Element Method (FEM). Both simulations show that the voltage and current behave as a propagating sinusoidal wave, which reflects when reaching the (right) end of the transmission line. However, in the FEM there is an overshoot when the wave reaches the end of the transmission line which causes peaks in the Hamiltonian (black line in Figure \ref{fig:ApproxH}), indicating that the power balance is not preserved in the FEM. Both methods show some high frequent oscillations, where the amplitude of the oscillations is approximately $2.3$ times higher for the method proposed here than for the FEM.  This is probably caused by numerical errors. As expected, for both models an increasing number of line segments yields smoother wave patterns.

\section{Conclusions}
\label{sec:conclusions}
This paper proposes a general method for approximating a port-Hamiltonian distributed parameter system on a 1D geometry by an aggregated interconnection of $N$ finite dimensional port-Hamiltonian systems. It is shown that an aggregated approximate model can be constructed that (i) inherits the Dirac structure of the distributed parameter model after projecting effort and flow spaces to finite dimensional subspaces, (ii) admits a Hamiltonian that coincides with the restriction of the Hamiltonian of the distributed parameter model on a finite dimensional subspace of the state of the distributed model. The aggregate model is port-Hamiltonian, of order $4N$ and preserves local and global power balances. It is lossless, linear, nonlinear, whenever the distributed parameter model is lossless, linear, nonlinear. It maintains the physical relevance of the energy density of the distributed model, and its state representation is sparse.

The theoretical development is based on the abstract  differential geometric setting of flow and effort variables of a generic one-dimensional nonlinear distributed parameter model in port-Hamiltonian form. The introduction of a mid-point $m$ in every element of the partition is an enabler to provide sufficient degrees of freedom in the design to meet the requirements of the aggregate model. Without the mid-point, it is not possible to meet the two requirements simultaneously. It is conjectured that this idea allows extensions to larger dimensional manifolds. The proposed method leaves freedom in the choice of shape functions defined on the manifold. Constructive examples of suitable shape functions in a one-dimensional geometry are given in Appendix A. 
Since shape functions have compact support on all elements of the partitioning, the state representation of the aggregate model is sparse. A simulation example clearly illustrates the preservation of the power and energy balance in a lossless transmission line, the sparsity of the state space representation of the approximate model, and the exact preservation of energy despite the numerical difficulty in simulating lossless systems. 

The method presented here preserves Dirac structures after spatial discretization and leaves time as a continuous variable. Extensions to discretizing both space and time can be based on the same Dirac structure, but are left as an important challenge. Additional challenges include questions on convergence of the Hamiltonian function $H_N$ and the aggregated model $\Sigma_N$ as $N\rightarrow\infty$. 

\bibliographystyle{ieeetr}        
\bibliography{port-Hamiltonian_discretization}

\appendix
\section{Construction of shape functions}
\label{App:shapeFunctions1D}
\setcounter{equation}{0}

\begin{figure}
\centering
\includegraphics[width=\columnwidth]{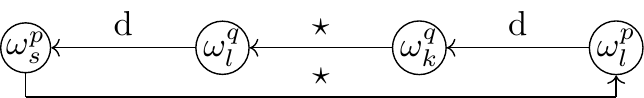}
\caption{Dependencies}
\label{fig:Dependencies1D}
\end{figure}

This Appendix covers an explicit construction of the different shape functions that are defined in the expansions \eqref{eq:fApprox1D} and \eqref{eq:eApprox1D}. The dependencies among 
the shape functions are derived in the same way as explained in \cite{Golo2004} and \cite{Tiemersma2016} and illustrated in Figure~\ref{fig:Dependencies1D} where the indices $k\in\{am,mb\}$, $s\in\{am,mb\}$ and $l\in\{a,m,b\}$.

Similar as in \cite{Golo2004}, the functions $\omega_{s}^p(z)$ for  $s \in \{am, mb\}$ are chosen as piecewise constant functions. Figure~\ref{fig:Dependencies1D} then implies that $\omega_l^q(z)$ with $l \in \{a,m,b\}$ are affine zero-forms. Subsequently, the Hodge star operator implies that $\omega_k^q(z)$ for $k \in \{am, mb\}$ are one-forms with the same affine structure as $\omega_l^q(z)$. Finally, the top right arrow in Figure~\ref{fig:Dependencies1D} implies that the zero-forms $\omega_l^p(z)$ with $l \in \{a,m,b\}$ become second order polynomials. The bottom arrow is not satisfied by these choices. 

An explicit construction of the shape functions is as follows.

\textbf{Definition of $\omega_{s}^p(z)$} \\
The shape functions $\omega_{s}^p(z)\in \Lambda^1(Z)$ indexed by $s \in \{am, mb\}$ are normalized according to \eqref{eq:1form1} and chosen piecewise constant one-forms according to
\begin{equation}
\omega_{s}^p(z) = \begin{cases} 
\begin{array}{lll}
\frac{\text{d}z_1}{\int_{Z_s} \text{d}z_1} &\text{ for } z \in Z_s \\
0 &\text{ for } z \notin Z_s,
\end{array}
\end{cases}
\label{eq:omega_s1D}
\end{equation} where $s \in \{am,mb\}$.
Let $m:=\frac{1}{2}(a+b)$ be the middle point of the element $Z_{ab}$. Then by \eqref{eq:omega_s1D}, $\omega_s^p(z)$ can be written 
as the one-form
\begin{equation*}
\omega_{s}^p(z) = \begin{cases} 
\begin{array}{lll}
\frac{2}{b-a} \text{ d}z_1 &\text{ for } z \in Z_s \\
0 &\text{ for } z \notin Z_s,
\end{array}
\end{cases}
\end{equation*} where $s \in \{am,mb\}$.

\textbf{Definition of $\omega_l^q(z)$} \\
By Figure~\ref{fig:Dependencies1D}, the functions $\omega_l^q(z)$  indexed by $l \in \{a,m,b\}$ then become affine zero-forms and satisfy the normalization conditions \eqref{eq:0form}. A general affine zero-form reads as $f(z) = \alpha z+\beta$ 
in the two real-valued parameters $\alpha$ and $\beta$. With the two normalization requirements on $\omega_a^q$ and $\omega_b^q$ this defines $\alpha$ and $\beta$ and sets
\begin{equation}
    \omega_a^q(z) = -\frac{1}{b-a}(z-b), \quad \omega_b^q(z) = \frac{1}{b-a}(z-a).
    \label{eq:omegaaq_omegabq}
\end{equation}
The function $\omega_m^q$ satisfies three requirements and is chosen to be a non-smooth, non-differentiable zero-form, which is affine on the interval $Z_{am}$ and affine on the interval $Z_{mb}$. It can be defined as
\begin{equation*}
    \omega_m^q(z) = \begin{cases}
    -\frac{2}{a-b}z+\frac{2a}{a-b} & \textrm{ for } z \in Z_{am} \\
    \frac{2}{a-b}z-\frac{2a}{a-b} & \textrm{ for } z \in Z_{mb}. \\
    \end{cases}
\end{equation*}

\textbf{Definition of $\omega_s^q(z)$} \\
Again by Figure~\ref{fig:Dependencies1D}, $\omega_s^q(z)$ with $s \in \{am, mb\}$ are chosen as affine one-forms that satisfy the following four requirements
\[
\int_{Z_{s_1}} \omega_{s_2}^q = \begin{cases}
1 \text{ for } s_1 = s_2 \\
0 \text{ for } s_1 \neq s_2,
\end{cases}
\]
where $s_1, s_2 \in \{am, mb\}$. Since a general affine one-form can be written as $f(z) = \alpha z + \beta \text{ d}z_1$,
there are two parameters to be determined from these two requirements. We infer
\begin{align*}
    &\omega_{am}^q(z) = \frac{1}{(a-b)^2}(-4z+a+3b) \textrm{d}z_1, \\ 
     &\omega_{mb}^q(z) = \frac{1}{(a-b)^2}(4z-3a-b) \textrm{d}z_1. 
\end{align*}

\begin{figure}
    \centering
    \includegraphics{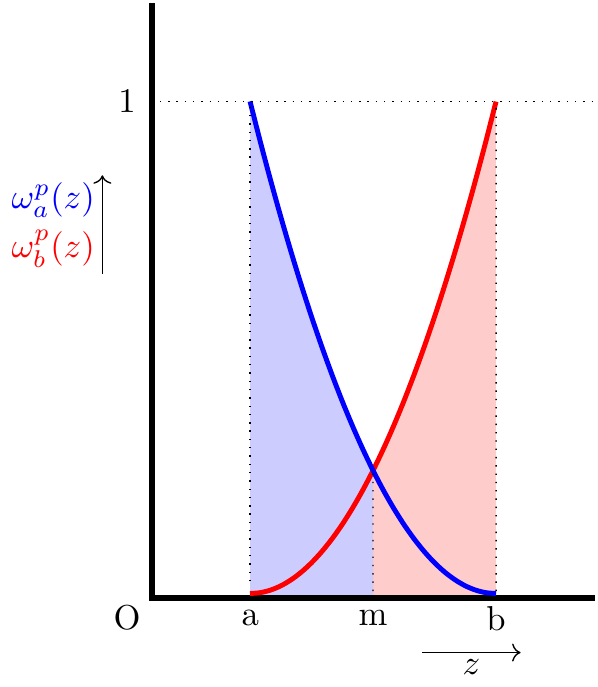}
    \caption{Area under $\omega_a^p(z)$ and under $\omega_b^p(z)$ are  equal. Here $\int_{Z_{am}} \omega_a^p \text{d}z$ is given in light blue and $\int_{Z_{mb}} \omega_b^p \text{d}z$ is in light red}
    \label{fig:sF}
\end{figure}

\textbf{Definition of $\omega_l^p(z)$} \\
By Figure~\ref{fig:Dependencies1D}, $\omega_l^p(z)$ for indices $l \in \{a,m,b\}$ become second order polynomials that satisfy  
\begin{equation*}
\omega^p_{l_1}(l_2) = \begin{cases}
1 \text{ for } l_1 = l_2 \\
0 \text{ for } l_1 \neq l_2, 
\end{cases}
\end{equation*}
where $l_1, l_2 \in \{a,b\}$ and
\begin{equation*}
\omega^p_{m}(l) = \begin{cases}
1 \text{ for } l = m \\
0 \text{ for } l \neq 0, 
\end{cases}
\end{equation*}
where $l \in \{a,m,b\}$. These functions therefore assume the form 
$f(z) = \alpha z^2 + \beta z + \gamma$
in the three parameters $\alpha$, $\beta$ and $\gamma$.  Since both $\omega_a^q$ and $\omega_b^q$ satisfy two normalization requirements, this is an underdetermined set of conditions. The 
condition
\[
    \omega_a^p(m) = \omega_b^p(m), \text{ and }
    0 < \omega_a^p(m) < 1
\]
imposes additional symmetry in the sense that $\omega_b^p$ is the same as $\omega_a^p$ mirrored at the line $z=m$ while
the integral of $\omega_a^p$ on $[a,m]$ equals the integral of $\omega_b^p$ on the interval $[m,b]$. See Figure~\ref{fig:sF}. 

Solving this by setting $\omega_a^p(m)=\frac{1}{4}$ this gives
\begin{align*}
    \omega_a^p(z) &= \frac{1}{(a-b)^2}\left( z^2-2bz+b^2 \right)\\
    \omega_b^p(z) &= \frac{1}{(a-b)^2} \left( z^2-2az+a^2 \right)\\
    \omega_m^p(z) &= \frac{1}{(a-b)^2} \left( -4z^2+4(a+b)z-4ab \right).
\end{align*}

\section{Proof of Proposition~\ref{prop:fdDirac}}
\label{App:Dirac1D}
\setcounter{equation}{0}

By \cite{Duindam2014}, to prove Proposition~\ref{prop:fdDirac} it suffices to show that
\begin{enumerate}
\item $FE^{\top}+EF^{\top} = 0$ and
\item rank$\begin{bmatrix}
F \mid E
\end{bmatrix} = $\textrm{dim} $\mathcal{F}$.
\end{enumerate}
In the construction of Proposition~\ref{prop:fdDirac}, condition (2) is equivalent to $\textrm{rank} \begin{bmatrix} F \mid E
\end{bmatrix} = \textrm{dim}(\mathcal{F})=6$ which is met. 

To verify condition (1), recall that $FE^{\top}+EF^{\top}$ can be written as
\begin{align*}
& FE^{\top}+EF^{\top} = \\
&=\begin{bmatrix}
M_e^{\top} & M_u^{\top}
\end{bmatrix} \begin{bmatrix}
M_f \\
M_y
\end{bmatrix} + \begin{bmatrix}
M_f^{\top} & M_y^{\top} 
\end{bmatrix} \begin{bmatrix}
M_e \\
M_u
\end{bmatrix} =\\
&= M_e^{\top}M_f + M_u^{\top}M_y + M_f^{\top}M_e+M_y^{\top}M_u.
\end{align*}
Using the expressions \eqref{eq:Mfdef},\eqref{eq:Medef} and substituting \eqref{eq:M21D}, \eqref{eq:M31D}, \eqref{eq:M41D}, \eqref{eq:M51D}, \eqref{eq:M61D}, \eqref{eq:Muydef} while setting $\sigma=0$ (the system is conservative) yields that $M_1=0$ and 
\begin{align*}
&FE^{\top}+EF^{\top} = \\
&\begin{bmatrix}
0 & M_4^\top M_2+M_3^\top M_5+M_6^\top \begin{bsmallmatrix}
-1 & 0 & 0\\
0 & 0 & 0 \\
0 & 0 & -1
\end{bsmallmatrix} \\
* & 0
\end{bmatrix}, 
\end{align*} where the star * denotes the transpose of the top right entry.
It thus suffices to show that \newline 
\mbox{$M_4^\top M_2+M_3^\top M_5+M_6^\top \begin{bsmallmatrix}
-1 & 0 & 0\\
0 & 0 & 0 \\
0 & 0 & -1
\end{bsmallmatrix}=0$}.
With $M_6$ as given in \eqref{eq:M61D}, it follows that
\begin{align*}
    &M_6^\top\begin{bmatrix}
-1 & 0 & 0\\
0 & 0 & 0 \\
0 & 0 & -1
\end{bmatrix} = \\
    &\int_{\partial Z_{ab}} \begin{bmatrix}
-\omega_a^q\mid_{\partial Z} \wedge \omega_a^p\mid_{\partial Z} & 0 & -\omega_b^q\mid_{\partial Z} \wedge \omega_a^p\mid_{\partial Z} \\
-\omega_a^q\mid_{\partial Z} \wedge \omega_m^p\mid_{\partial Z} & 0 & 
-\omega_b^q\mid_{\partial Z} \wedge \omega_m^p\mid_{\partial Z} \\
-\omega_a^q\mid_{\partial Z} \wedge \omega_b^p\mid_{\partial Z} & 0 &
-\omega_b^q\mid_{\partial Z} \wedge \omega_b^p\mid_{\partial Z} \\
\end{bmatrix}
\end{align*}
Since these shape functions satisfy \eqref{eq:0form}, we can conclude this is equivalent to 
\[
    M_6^\top\begin{bmatrix}
-1 & 0 & 0\\
0 & 0 & 0 \\
0 & 0 & -1
\end{bmatrix} = \begin{bmatrix}
1 & 0 & 0 \\
0 & 0 & 0 \\
0 & 0 & -1 \\
\end{bmatrix}.
\]

Next, the expressions for the shape functions derived in Appendix~\ref{App:shapeFunctions1D} are substituted in the matrices $M_2, M_3, M_4, M_5$. An analytic, symbolic or numerical computation then reveals that
\begin{equation*}
    M_4^\top M_2+M_3^\top M_5= \begin{bmatrix}
    -1 & 0 & 0 \\
    0 & 0 & 0 \\
    0 & 0 & 1
    \end{bmatrix}
\end{equation*}
showing that $FE^{\top}+EF^{\top} =0$. 

\end{document}